\documentclass[english]{cccconf}
\usepackage[comma,numbers,square,sort&compress]{natbib}

\usepackage{color, colortbl}
\definecolor{Gray}{gray}{0.9}

\usepackage{stfloats}

\usepackage{cuted}

\usepackage{mathrsfs}
\usepackage{amsthm}

\newtheorem{lemma}{Lemma}

\usepackage{hyperref} 

\usepackage{algorithm,algorithmic}

\usepackage{amsfonts}
\usepackage{caption}
\usepackage{graphicx}

\usepackage[utf8]{inputenc}
\usepackage{amssymb}

\usepackage{nomencl} 
\makenomenclature
\usepackage{etoolbox}
\renewcommand{\nomgroup}[1]{%
  \item[\bfseries
  \ifstrequal{#1}{A}{Given Parameters}{%
  \ifstrequal{#1}{B}{Decision variables}{%
  \ifstrequal{#1}{C}{Other Symbols}{}}}%
]}
\usepackage{amsmath}
\usepackage{algorithmic}

\usepackage{array}

\usepackage{subfigure}

\usepackage{float}

\usepackage{stfloats}

\usepackage{nomencl} 
\makenomenclature
\usepackage{etoolbox}
\renewcommand{\nomgroup}[1]{%
	\item[\bfseries
	\ifstrequal{#1}{A}{Given Parameters}{%
		\ifstrequal{#1}{B}{Decision variables}{%
			\ifstrequal{#1}{C}{Other Symbols}{}}}%
	]}

\usepackage{mathtools}
\DeclarePairedDelimiter\ceil{\lceil}{\rceil}

\begin{document}

\title{Evacuation Problem Under the Nuclear Leakage Accident}

\author{Canqi Yao\aref{hit,sustech},
        Shibo Chen\aref{sustech},
        Zaiyue Yang\aref{sustech}}


\affiliation[hit]{School of Mechatronics Engineering, Harbin Institute of Technology, Harbin, 150000, China
        \email{vulcanyao@gmail.com}}
\affiliation[sustech]{Shenzhen Key Laboratory of Biomimetic Robotics and Intelligent Systems, Department of Mechanical and Energy Engineering, and the Guangdong Provincial Key Laboratory of Human-Augmentation and Rehabilitation Robotics in Universities, Southern University of Science and Technology, Shenzhen 518055, China
        \email{chensb@sustech.edu.cn, yangzy3@sustech.edu.cn}}

\maketitle


\begin{abstract}
To handle the detrimental effects brought by leakage of radioactive gases at nuclear power station, we propose a bus based evacuation optimization problem. The proposed model incorporates the following four constraints, 1) the maximum dose of radiation per evacuee, 2) the limitation of bus capacity, 3) the number of evacuees at demand node (bus pickup stop), 4) evacuees balance at demand and shelter nodes, which is formulated as a mixed integer nonlinear programming (MINLP) problem. Then, to eliminate the difficulties of choosing a proper M value in Big-M method, a Big-M free method is employed to linearize the nonlinear terms of the MINLP problem. Finally, the resultant mixed integer linear program (MILP) problem is solvable with efficient commercial solvers such as CPLEX or Gurobi, which guarantees the optimal evacuation plan obtained. To evaluate the effectiveness of proposed evacuation model, we test our model on two different scenarios (a random one and a practical scenario). For both scenarios, our model attains executable evacuation plan within given 3600 seconds computation time.  \end{abstract}


\footnotetext{ Accepted by 2021 40th Chinese Control Conference (CCC). IEEE, 2021. This work was supported in part by the National Key Research and
Development Program of China under Grant 2019YFB1705401; in part by
the Natural Science Foundation of China under Grant 61873118 and Grant
61903179; in part by the Science, Technology and Innovation Commission
of Shenzhen Municipality under Grant ZDSYS20200811143601004 and Grant
RCBS20200714114918137}

\makenomenclature
\mbox{}

\nomenclature[A]{$G(N,A)$}{Directed graph, \textit{N} denotes the nodes in the digraph and \textit{A} denotes arcs in the digraph connecting nodes.}
\nomenclature[A]{$\mathcal{T}_{ij}$}{Travel time between node \textit{i} and \textit{j},$(i,j)\in A$}

\nomenclature[A]{$\tau_{ij}$}{ Nuclear radiation per second ($mSv/s$) between node $i,j$, $(i,j)\in A$  }
\nomenclature[A]{$\eta_i$}{Nuclear radiation per second ($mSv/s$) at node $i$}
\nomenclature[A]{$\sigma$}{The maximum dose of nuclear radiation for evacuees}
\nomenclature[A]{Q}{The capacity of bus}
\nomenclature[A]{$D_j$}{Number of evacuees at demand node $j, j\in P$}

\nomenclature[B]{$x_{ij}^{mt}$}{Binary decision variable that equals  1 if trip \textit{t} for bus \textit{m} traverses arc $(i,j)$, else 0, $\forall (i,j)\in A,m\in V, \quad  t\in \mathscr{T}$ }
\nomenclature[B]{$\mathcal{T}^{mt}_{i}$}{Visiting time for the $m^{th}$ bus arriving node \textit{i} at the $t^{th}$ trip, $\forall i\in N,m\in V,t\in \mathscr{T}$}
\nomenclature[B]{$b_i^{mt}$}{Number of evacuees from \textit{i} assigned to (or, if \textit{i} is a shelter, released from) bus \textit{m} after trip $t$, $\forall i\in P\cup S,m\in V,t\in \mathscr{T}$}


\printnomenclature

\section{Introduction}

\par According to International atomic energy agency (IAEA), a nuclear and radiation accident could lead to significant consequences to people, and the environment. The most severe nuclear disasters, such as the Chernobyl disaster in 1986\cite{hobeika1994decision} and Fukushima Daiichi nuclear disaster\cite{hindmarsh2013nuclear} in 2011, have caused enormous loss in human life and cost with lasting release of radioactive gases into the surrounding environment. To cope with these terrible accidents, there are some protective methods identified to substantially reduce or avert exposures of radioactive gases. Among these methods, evacuation, the frequently used protective method in nuclear disaster, is to transport evacuees from the affected area to avoid receiving too much radiation.\par Motivated by the challenging emergency problems in natural or man-created disasters like hurricanes, earthquakes, flooding, landslide, and nuclear accidents, evacuation planning has been investigated extensively both in academic and industry\cite{oh2019efficient,goerigk2013branch,talebi1985stochastic}.  Oh et al.\cite{oh2019efficient} propose a new network considering congestion called a multi-class time-expanded (MCTE) network under which an efficient heuristic algorithm is devised to solve the evacuation problem.  In the framework of branch and bound method, goerigk et al. \cite{goerigk2013branch} propose 
multiple approaches for finding both lower and upper bounds for the bus evacuation problem with several node pruning techniques and branching rules. Dhamala gives a systematic analysis of network flow models used in emergency evacuation and their applications in a survey paper \cite{dhamala2015survey}. For more surveys concerning with evacuation planning, readers refer to \cite{chalmet1982network,yusoff2008optimization,hamacher2002mathematical,altay2006or,schadschneider2008evacuation}.\par Different from the previously mentioned evacuation problem, under the circumstance of nuclear leakage, the nuclear radiation has to be considered into the evacuation problem to make sure the safety of every evacuee. As for the optimization of emergency evacuation from
nuclear accidents, zou et al.\cite{zou2018optimization}  propose a solution method for solving the shortest evacuation path problem based on
graph theory. However, this solution method has limited applications and can not be applied into complicated scenarios. To handle the uncertainties such as fuzzy, interval and fuzzy random variables in an evacuation management system, huang et al. \cite{huang2017inexact} devise an inexact fuzzy stochastic
chance constrained programming (IFSCCP) method.\par In this paper, a bus based evacuation optimization problem is proposed incorporating the following constraints: maximum radiation value of the evacuees, the limitation of bus capacity, the number of evacuees at demand node ( bus pickup stop ), evacuees balance at demand nodes and shelter nodes, which is formulated as a mixed integer nonlinear programming (MINLP) problem. Then, instead of adopting Big-M method which has difficulties in choosing a proper M value, a Big-M free method is proposed to linearize  the nonlinear terms of the MINLP problem. Finally, the resultant mixed integer linear program (MILP) problem is solved by efficient commercial solvers CPLEX to obtain the optimal evacuation plan.  

\par The rest of the paper is organized as follows. In section II, we propose a bus based evacuation problem formulated as a mixed integer bi-linear program (MIBP) problem. To cope with the non-linearity in MIBP problem, a linearization method is proposed in section III. Numerical experiments are conducted in section IV. We conclude our paper in section V.

\section{Bus based evacuation problem }

In this section, we consider a directed graph  $G(N, A)$, where $N$ and $A$ denote the set of nodes and arcs respectively. Specifically, $N=\{D,P,S\}$ is composed of three subsets of nodes: $D$, a set of depot nodes where buses are initially located and dispatched from; $P$, a set of demand nodes, each of which represents a pickup location serving a neighborhood or facility requiring evacuation services; and $S$, a set of shelter nodes. $A=\{(i,j)|i\in D, j\in P\}\cup\{(i,j)|i\in P, j\in S\}\cup\{(i,j)|i\in S, j\in P\}$ specifies the arcs connecting nodes in $N$.
We keep track of consecutive trips of buses  with  index
 $\mathscr{T}=\{1,\ldots,T\}$.
 
\par
\begin{figure}[ht]
\centering
\includegraphics[width=.95\linewidth]{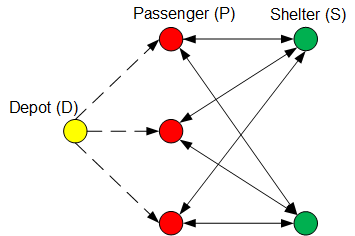}
\caption{Evacuation Network}
\label{ENetwork}
\end{figure}

 A toy example about the evacuation network is shown in Fig.\ref{ENetwork}. The evacuation buses start from depot (yellow node), pick the evacuees at passenger nodes (red nodes) and deliver these evacuees to shelters (green nodes).  Then, evacuation buses return to passenger nodes to pick 
the remaining evacuees until all evacuees
are rescued.

\subsection{Mathematical Formulation }
In this subsection, the objective function and constraints related to evacuation problem under nuclear leakage are clarified and explained in details. The objective function focuses on the minimization of travel time of all buses over all trips:
\begin{equation}\notag
\begin{aligned}
 \min_{x_{ij}^{mt},\mathcal{T}_i^{mt},b_i^{mt}}  \sum_{(i,j)\in A}\sum_{m\in V}  \sum_{t\in\mathscr{T}} \mathcal{T}_{ij}  x_{ij}^{mt}  
\end{aligned}
\end{equation}



\par  The total received radiation of evacuees consisting of two types of radiation, waiting radiation and escape radiation, is limited under a given upper bound ($\delta$)  to ensure the safety of evacuees in (\ref{Radiation}).
\begin{equation}\label{Radiation}
\sum_{t\in\mathscr{T}}\sum_{ j \in S}  \tau_{i j} \mathcal{T}_{ij}  b_i^{mt}  +  \sum_{t\in\mathscr{T}} \eta_i \mathcal{T}^{mt}_i b_i^{mt} \leq  \delta  , \quad \forall i\in P, m \in V
\end{equation}



 Constraint (\ref{time}) specifies that the evacuation plan should satisfy the time constraints.
\begin{equation}\label{time}
\begin{aligned}
    \mathcal{T}^{mt+1}_j &\geq \mathcal{T}^{mt}_i+x_{ij}^{mt}\mathcal{T}_{ij},\\& \forall i\in N, j\in N, m\in V,t\in\mathscr{T}\setminus T
\end{aligned}
\end{equation}


Constraint (\ref{fbc_demand}) is the flow-balance constraint for the demand nodes; it ensures that a bus travels to demand node $j$ on trip \textit{t} leaves node $j$ on trip \textit{t + 1};
Constraint (\ref{fbc_shelter}) is the flow-balance constraint for the shelters; it does not require the bus to leave the shelter \textit{i}, that is, the last trip a bus makes can end at a shelter. Constraint (\ref{c_served}) allows a bus to make at most one trip at a time.

\begin{equation}\label{fbc_demand}
\begin{aligned}
\sum_{i:(i, j) \in A} x_{i j}^{m t}&=\sum_{k:(j, k) \in A} x_{j k}^{m(t+1)},\\ &\quad \forall j \in P, m \in V, t\in\mathscr{T}\setminus T 
\end{aligned}
\end{equation}

\begin{equation}\label{fbc_shelter}
\begin{aligned}
 \sum_{i:(i, j) \in A} x_{i j}^{m t} &\geq \sum_{k:(j, k) \in A} x_{j k}^{m(t+1)},\\& \quad \forall j \in S, m \in V, t\in\mathscr{T}\setminus T  
\end{aligned}
\end{equation}

\begin{equation}\label{c_served}
\sum_{(i, j) \in A} x_{i j}^{m t} \leq 1, \quad \forall m \in V, t\in\mathscr{T}
\end{equation}


 Constraint (\ref{c_1leave}) specifies that the first trip of each bus start from its depot, and Constraint (\ref{c_laterleave}) ensures that the buses do not leave the depot for later trips. Constraint (\ref{c_notend}) does not allow the last trip a bus can make to end at a demand node. 
 
 \begin{equation}\label{c_1leave}
\sum_{j\in N}x_{i j}^{m 1}=1, \quad \forall i \in Y, m \in V_{i}
\end{equation}

\begin{equation}\label{c_laterleave}
x_{i j}^{m t}=0, \quad \forall i \in Y, j \in N, m \in V, t\in\mathscr{T}\setminus 1
\end{equation}

\begin{equation}\label{c_notend}
x_{i j}^{m T}=0, \quad \forall j \in P, i \in S, m \in V
\end{equation}

 Constraint (\ref{c_capacity1}) dictates that a bus can only pick up evacuees from node $j$ if it has in fact traveled to node $j$. Constraint (\ref{c_capacity2}) is the bus capacity constraint.  Eq. (\ref{c_demand}) ensures that all evacuees are picked up.
 
 \begin{equation}\label{c_capacity1}
b_{j}^{m t} \leq \sum_{(i, j) \in A} Q x_{i j}^{m t}, \quad \forall j \in N, m \in V, t\in\mathscr{T}
\end{equation}

\begin{equation}\label{c_capacity2}
0 \leq \sum_{j \in P} \sum_{l=1}^{t} b_{j}^{m l}-\sum_{k \in S} \sum_{l=1}^{t} b_{k}^{m l} \leq Q, \quad \forall m \in V, t\in\mathscr{T}
\end{equation}


\begin{equation}\label{c_demand}
\sum_{m \in V} \sum_{t=1}^{T} b_{j}^{m t}=D_{j}, \quad \forall j \in P
\end{equation}

 Constraint (\ref{c_shelterDemand}) ensures that all evacuees have to be  delivered to shelters.

\begin{equation}\label{c_shelterDemand}
\sum_{j \in P} \sum_{t=1}^{T} b_{j}^{m t}=\sum_{k \in S} \sum_{t=1}^{T} b_{k}^{m t}, \quad \forall m \in V
\end{equation}

Constraints (\ref{c_x}) and (\ref{c_b}) are the logical binary and integer restrictions on $x$ and $b$, respectively.  Constraint (\ref{c_T}) specifies the non-negativity of $\mathcal{T}$.

\begin{equation}\label{c_x}
x_{i j}^{m t} \in\{0,1\}, \quad \forall(i, j) \in A, m \in V, t\in\mathscr{T}
\end{equation}

\begin{equation}\label{c_b}
b_{j}^{m t}\in \mathbb{Z} \geq 0, ,\quad \forall j \in N, m \in V, t\in\mathscr{T}
\end{equation}

\begin{equation}\label{c_T}
\mathcal{T}_{j}^{m t} \geq 0, \quad \forall j \in N, m \in V, t\in\mathscr{T}
\end{equation}

Combining all constraints and objective functions, the evacuation problem is mathematically formulated as problem 1, a mixed integer bi-linear programming problem. Due to the difficulties of these bi-linear terms $\mathcal{T}_i^{mt} b_i^{mt}$ in (2), commercial solvers such as Gurobi or CPLEX  cannot be directly employed to solve problem 1. 

\par \textbf{Problem 1 (MIBP)}
\begin{equation}\notag
\begin{aligned}
\min_{x_{ij}^{mt},\mathcal{T}_i^{mt},b_i^{mt}}   \sum_{(i,j)\in A}\sum_{m\in V}  \sum_{t\in\mathscr{T}} \mathcal{T}_{ij}  x_{ij}^{mt}  
\end{aligned}
\end{equation}
$$\textit{s.t.}\qquad \text{(\ref{Radiation})-(\ref{c_T})}$$

In the next section, to solve the bus evacuation problem,we exploit an exact linearization method to transform problem 1 as the mixed integer linear programming problem (MILP) equivalently.

\section{MILP Reformulation}
In this section, in order to deal with the integer bi-linear terms in constraints, inspired by \cite{Pereira2005}, we encode the integer variable $b_i^{mt}$ with a set of binary variables. Then, the integer bi-linear terms are represented by a set of binary bi-linear terms, which can be  exactly linearized. As for the resultant binary bi-linear terms, instead of employing traditional Big-M method to linearize binary bi-linear terms\cite{cerna2017optimal}, a novel Big-M free method is devised, thus avoids the difficulty in choosing an appropriate M value.

\subsection{Binary Expansion Method}
The basic idea lies in expressing the integer number in binary
numeral system. For example, suppose that the plausible values for the number of evacuees $b_i^{mt}$ are in the range  $[0,Q]$. The
integer number is expressed as
\begin{equation*}
    b_i^{mt}=\sum_{n\in \mathscr{N}} 2^n y_n^{imt}
\end{equation*}
where $\mathscr{N}=\{1,\ldots,num\}$, $num=\ceil{log_2(1+Q)-1}$ and $y_n^{imt}\in\{0,1\}$.

With the binary-represented integer number, we can reformulate the integer bi-linear terms $\mathcal{T}_i^{mt} b_i^{mt}$ in (\ref{Radiation}) as a set of binary bi-linear terms  
\begin{equation*}
 \mathcal{T}_i^{mt} b_i^{mt}= \mathcal{T}_i^{mt} \sum_{n\in \mathscr{N}} 2^n y_n^{imt}= \sum_{n\in \mathscr{N}} 2^n \mathcal{T}_i^{mt}  y_n^{imt} . 
\end{equation*}
Finally, (\ref{Radiation}) is further reduced to the following constraint. 
\begin{equation}\label{Radiation_B}
\begin{aligned}
\sum_{t\in\mathscr{T}}\sum_{ j \in S}  \tau_{i j} \mathcal{T}_{ij}  b_i^{mt}  + & \sum_{t\in\mathscr{T}} \eta_i \sum_{n\in \mathscr{N}} 2^n \mathcal{T}_i^{mt}  y_n^{imt} \leq  \delta  , \\&\forall i\in P, m \in V
\end{aligned}
\end{equation}

\subsection{Big-M Free Linearization Method}
Recently, Big-M method prevails in the linearization method of bi-linear term \cite{chen2016optimal,cerna2017optimal,james2017autonomous}. {Note that, defining a proper \textit{M} is vital in the Big-M method and a large \textit{M} may result in serious numerical difficulties on a computer, while a small \textit{M} may not guarantee the optimal solution} \cite{pourakbari2015unambiguous,pineda2019solving}. Inspired from \cite{2012multilinear}, to overcome the inconveniences brought by Big-M method, we propose a Big-M free method without introducing a big \textit{M}.

\begin{lemma}\label{Equi}
Given an binary bi-linear term $\omega=x_1 x_2$, $x_1\in [x^L_1,x^U_1]$, $x_2\in \{0,1\}$, if $x^L_1=0$, $\omega$ can be exactly linearized by the following constraints.
\begin{subequations}\label{Bilinear}
\begin{align}
x_2-1\leq \omega\leq x_2 x^R_1\\    
x^L_1\leq \omega\leq x^R_1
\end{align}
\end{subequations}
\end{lemma}

\begin{proof} The equivalence between binary bi-linear term and (\ref{Bilinear}) is clarified as follows
\begin{itemize}
    \item If $x_2=0$, $\omega=0$ and  (\ref{Bilinear}) is reduced to $-1\leq \omega\leq 0,  
0\leq \omega\leq x^R_1$, implying $\omega=0$.

    \item If $x_2=1$, $\omega=x_1$ and  (\ref{Bilinear}) is reduced to $0\leq \omega\leq x^R_1,  
0\leq \omega\leq x^R_1$. 
\end{itemize}
Proof completes.
\end{proof}

Along with (\ref{M_Radiation})-(\ref{Bound}), $v^{imt}_n$ is introduced to eliminate these bi-linear terms in problem 1. Let


\begin{equation}\label{Binary}
v^{imt}_n=2^n \mathcal{T}_i^{mt}  y_n^{imt},\forall n\in \mathscr{N}; i\in N; m\in V; t\in \mathscr{T}
\end{equation}

\subsubsection{Linearization of  radiation constraints}
Following Lemma \ref{Equi},  (\ref{Radiation_B}) is equivalently transformed as (\ref{M_Radiation})-(\ref{Bound}).


\begin{equation}\label{M_Radiation}
 \sum_{t\in\mathscr{T}}\sum_{ j \in S}  \tau_{i j} \mathcal{T}_{ij}  b_i^{mt}  +  \sum_{t\in\mathscr{T}} \eta_i \sum_{n\in \mathscr{N}} v_n^{imt} \leq  \delta  , \quad \forall i\in P, m \in V
\end{equation}


\begin{equation}\label{Bound0}
y_n^{imt}-1 \leq  v_n^{imt} \leq y_n^{imt} 2^n\mathcal{\bar{T}}_i^{mt},   \forall n\in\mathscr{N}, i\in P, m\in V,\\ t\in\mathscr{T}  
\end{equation}

\begin{equation}\label{Bound}
0\leq  v_n^{imt} \leq 2^n\mathcal{\bar{T}}_i^{mt} ,   \forall n\in\mathscr{N}, i\in P, m\in V,\\ t\in\mathscr{T}  
\end{equation}



\par Finally, the time-consuming problem 1 (MIBP) can be exactly simplified as problem 2 (MILP). Unlike problem 1, problem 2 is a mixed integer linear program problem, which can be directly solved by commercial solvers. 

\textbf{Problem 2 (MILP)}
\begin{equation}\notag
\begin{aligned}
\min_{x_{ij}^{mt},y_n^{imt},\mathcal{T}_i^{mt},b_i^{mt}}    \sum_{(i,j)\in A}\sum_{m\in V}  \sum_{t\in\mathscr{T}} \mathcal{T}_{ij}  x_{ij}^{mt}  
\end{aligned}
\end{equation}
$$\textit{s.t.}\qquad \text{(\ref{time})-(\ref{c_T}) and} (\ref{Binary})-(\ref{Bound})$$


\section{Numerical Results}
In order to evaluate the effectiveness of proposed mathematical model of evacuation problem, extensive simulations are carried out. The proposed mathematical model is firstly tested on a random generated scenario. Then to test the performance of our model on practical case, we put our model on a practical scenario, which is generated based on the road map of a small town in Chile\footnote{\url{https://github.com/scespinoza/bus-evacuation-problem/tree/master/instances/paipote}}.  Besides, all scenarios are tested using the DOcplex API on a Python 3.7 environment with a time limit of 3600 seconds, using a PC with an Intel Core i5-7500 3.4GHz CPU and 8GB of RAM.
\subsection{Random Scenario}
\subsubsection{Parameter settings}
Proposed MILP formulation is tested on 7 instances (random generated) detailed in Table I. The differences between these instances lie in the number of shelters, the number of buses and capacity of buses.

\begin{table}[ht]
\label{table_inst1}
\begin{center}
\caption{Parameter Settings \textit{(Random Scenario)} }
\begin{tabular}{ ccccccc } 
\hline
\hline
Instance & Depot & Pickups & Shelters & Buses & Q\\ 
\hline
\hline
1& 1 & 5 & 2 & 20  & 25\\
\rowcolor{Gray}
2& 1 & 5 & 2 & 20  & 30\\
3& 1 & 5 & 2 & 25  & 30\\
\rowcolor{Gray}
4& 1 & 5 & 2 & 30  & 30\\
5& 1 & 5 & 3 & 20  & 30\\
\rowcolor{Gray}

6& 1 & 5 & 3 & 25  & 30\\
7& 1 & 5 & 3 & 30  & 30\\
\hline
\hline
\end{tabular}
\end{center}
\end{table}

\begin{table}[ht]
\begin{center}
\caption{Simulation Results \textit{(Random Scenario)} }
\begin{tabular}{ccccccc} 
\hline\hline
Instance & Optimality Gap & Elapsed Time (s) & $T_{evac}$ & Cost\\ 
\hline\hline

1&  0.01\% & 2380.11 & 1478.43  & 11746.74   \\
\rowcolor{Gray}
2&  0.01 \% & 1661.77 & 923.41 & 9808.52  \\
3&  0\% & 12.77 & 1204.71  & 10383.25   \\
\rowcolor{Gray}
4&  0\% & 94.53 & 923.41 & 11033.61   \\
5&   $\boldsymbol{2.25\% }$ & time limit & 853.71 & 9227.99  \\
\rowcolor{Gray}
6&  0\% & 91.78 & 923.41  & 9266.84   \\
7&  0\% & 74.42 & 923.41 & 10017.19   \\
\rowcolor{Gray}

Average&  0.32\% &  1130.76 & / & /   \\

\hline\hline
\end{tabular}
\end{center}
\label{table_res1l}
\end{table}

\subsubsection{Experimental results}

\begin{figure}[ht]
\centering
\includegraphics[width=.85\linewidth]{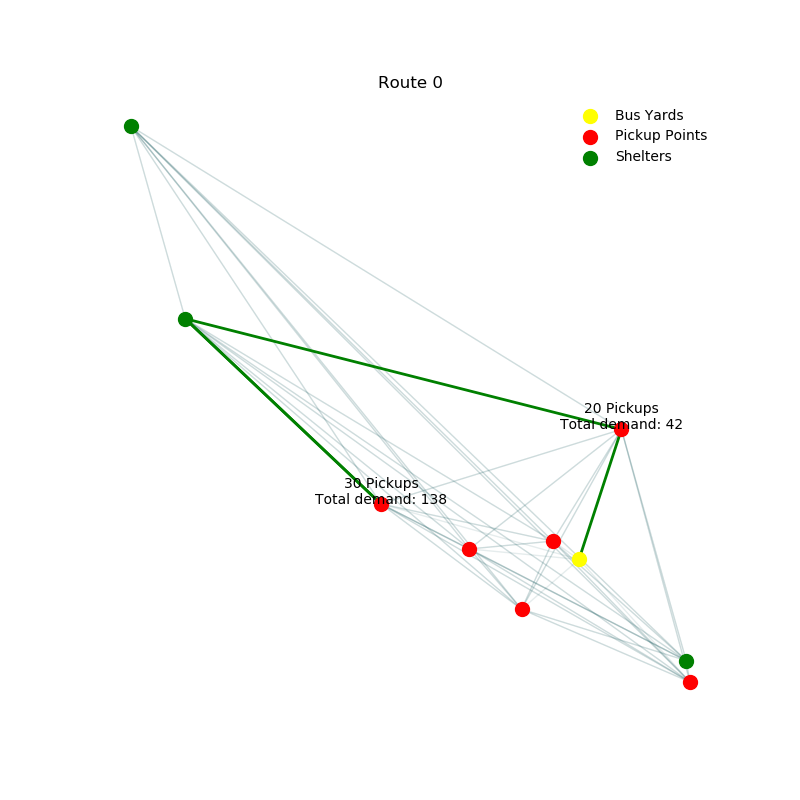}
\caption{Evacuation route of vehicle 1}
\label{R1}
\end{figure}

\begin{figure}[ht]
\centering
\includegraphics[width=.85\linewidth]{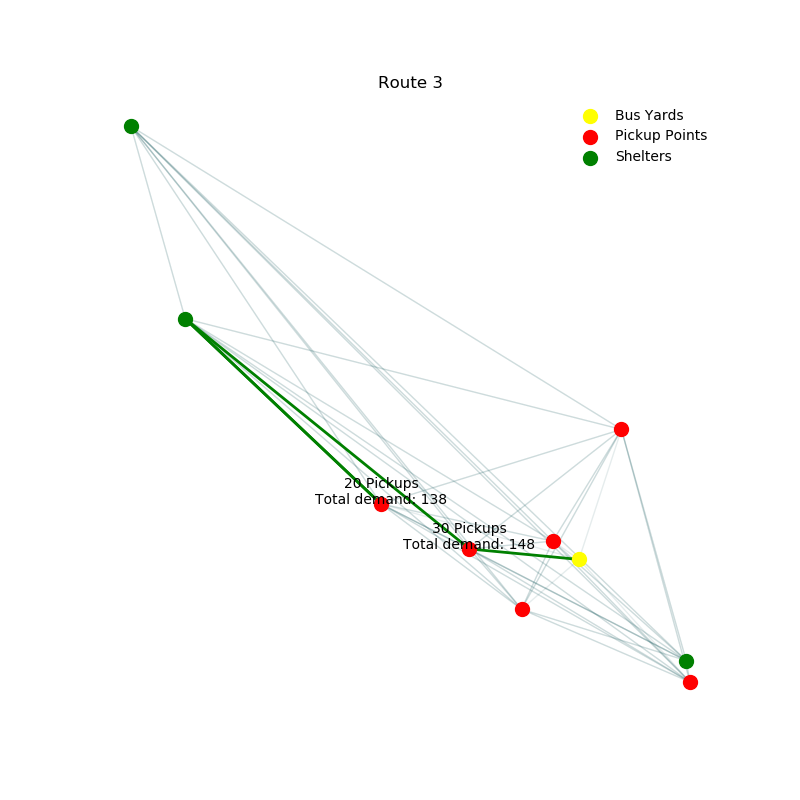}
\caption{Evacuation route of vehicle 4}
\label{R2}
\end{figure}

\begin{figure}[ht]
\centering
\includegraphics[width=.85\linewidth]{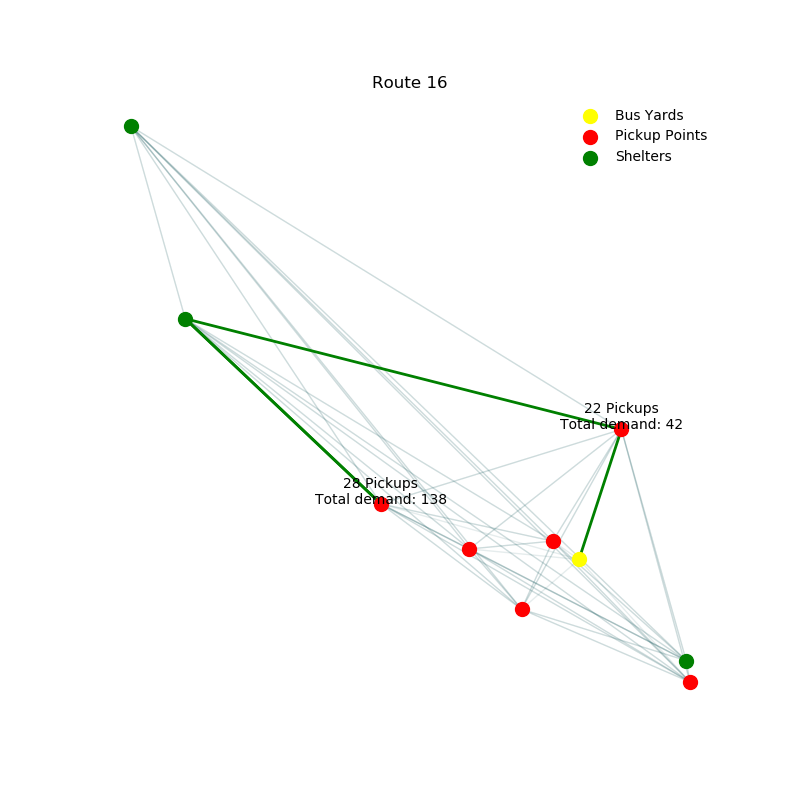}
\caption{Evacuation route of vehicle 17}
\label{R3}
\end{figure}

The simulation results of 7 instances are shown in Table II and compared in four categories, optimality gap, elapsed time, evacuation time\footnote{The longest evacuation time among all buses.}, objective value\footnote{The sum of evacuation time of all buses.}. Note that 6 out of 7 instances obtain an optimal solution within prescribed elapsed time in which 4 out 7 instances attain its optimal solution under 150 seconds, and the upper bound of elapsed time is set as 3600 seconds (1 hour). The average optimality loss for 7 instances is $0.32\% $, and the highest optimality gap is $2.25\% $ occurred in Instance 5. The detailed evacuation routes for vehicle 1, vehicle 4 ,and vehicle 17 of instance 2  
are illustrated in Fig.\ref{R1}, Fig.\ref{R2}, and Fig.\ref{R3} respectively.

\subsection{Practical Scenario}

\subsubsection{Parameter settings} The major aspects of practical scenario are the same with the random one, except that the differences exist in the structures of road map, the number of pickup nodes and shelters.

\begin{table}[ht]
\label{table_inst}
\begin{center}
\caption{Parameter Settings \textit{(Practical Scenario)} }
\begin{tabular}{ ccccccc } 
\hline
\hline
Instance & Depot & Pickups & Shelters & Buses & Q\\ 
\hline
\hline
1& 1 & 6 & 3 & 20  & 25\\
\rowcolor{Gray}
2& 1 & 6 & 3 & 20  & 30\\
3& 1 & 6 & 3 & 25  & 30\\
\rowcolor{Gray}
4& 1 & 6 & 3 & 30  & 30\\
5& 1 & 6 & 2 & 20  & 30\\
\rowcolor{Gray}
6& 1 & 6 & 2 & 25  & 30\\
7& 1 & 6 & 2 & 30  & 30\\
\hline
\hline
\end{tabular}
\end{center}
\end{table}

\subsubsection{Experimental results}

\begin{table}[ht]
\begin{center}
\caption{Simulation Results (\textit{Practical Scenario})}
\begin{tabular}{ccccccc} 
\hline\hline
Instance & Optimality Gap & Elapsed Time (s) & $T_{evac}$ & Cost\\ 
\hline\hline
1&  0.01\% & 1964.39 & 626.92  & 7362.08   \\
\rowcolor{Gray}
2&  $\boldsymbol{1.64 \%}$  & time limit & 654.03 & 5932.21  \\
3&  0\% & 0.39 &330.84  &5881.30   \\
\rowcolor{Gray}
4&  0\% & 0.56 & 330.85 &6827.18   \\
5&   0.01\% & 109.70 & 654.03 & 5932.21  \\
\rowcolor{Gray}
6&  0\% & 0.63 &330.84  &5881.35   \\
7&  0\% & 0.46 & 330.85 &6827.18   \\
\rowcolor{Gray}
Average&  0.24\% &  810.86 & / & /   \\

\hline\hline
\end{tabular}
\end{center}
\label{table_resl}
\end{table}

\begin{figure}[ht]
\centering
\includegraphics[width=.85\linewidth]{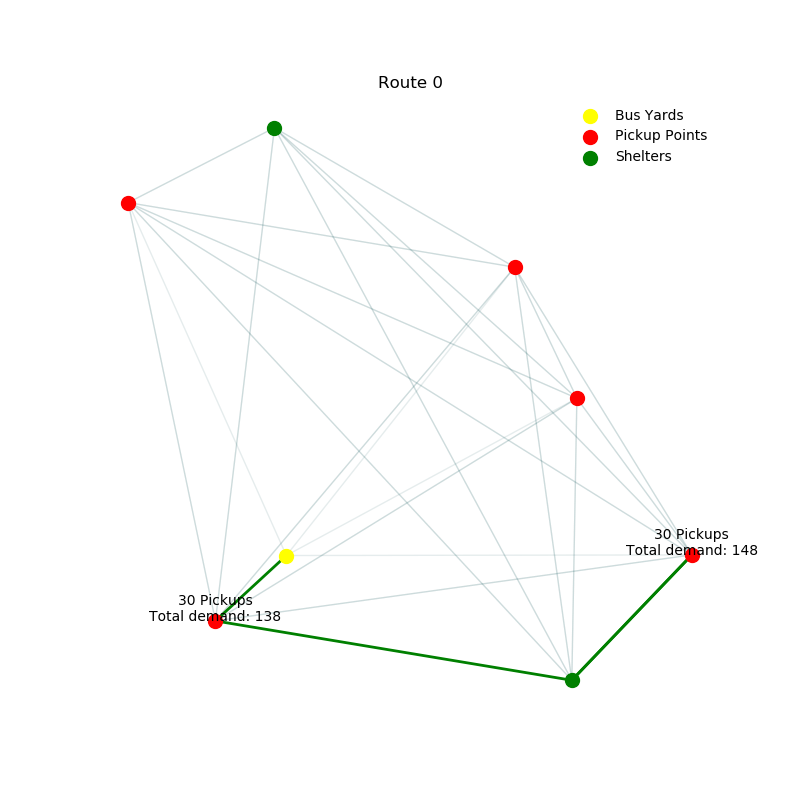}
\caption{Evacuation route of vehicle 1}
\label{P1}
\end{figure}

\begin{figure}[ht]
\centering
\includegraphics[width=.85\linewidth]{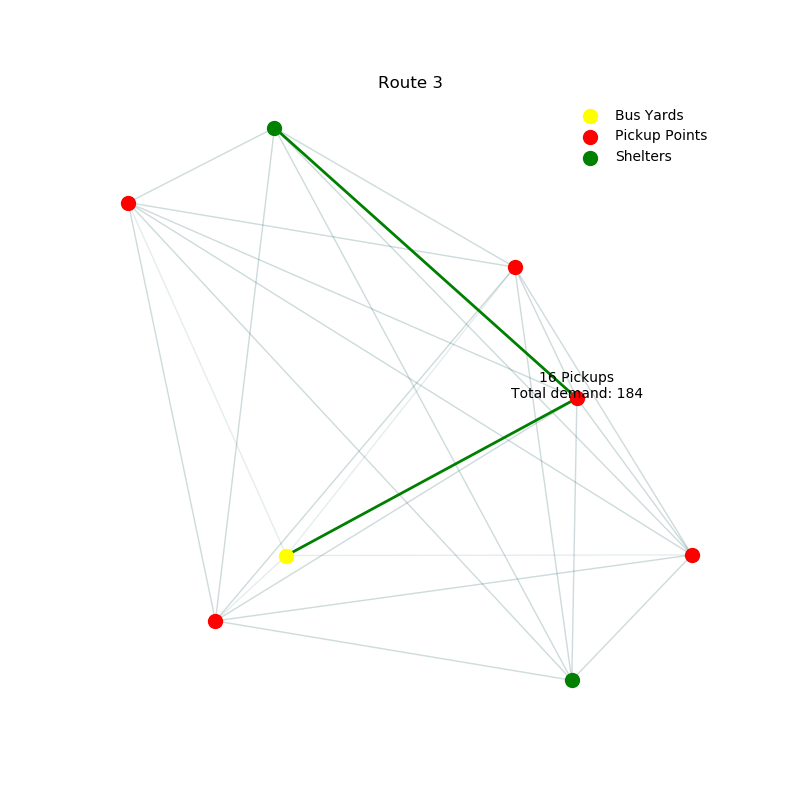}
\caption{Evacuation route of vehicle 4}
\label{P2}
\end{figure}

\begin{figure}[ht]
\centering
\includegraphics[width=.85\linewidth]{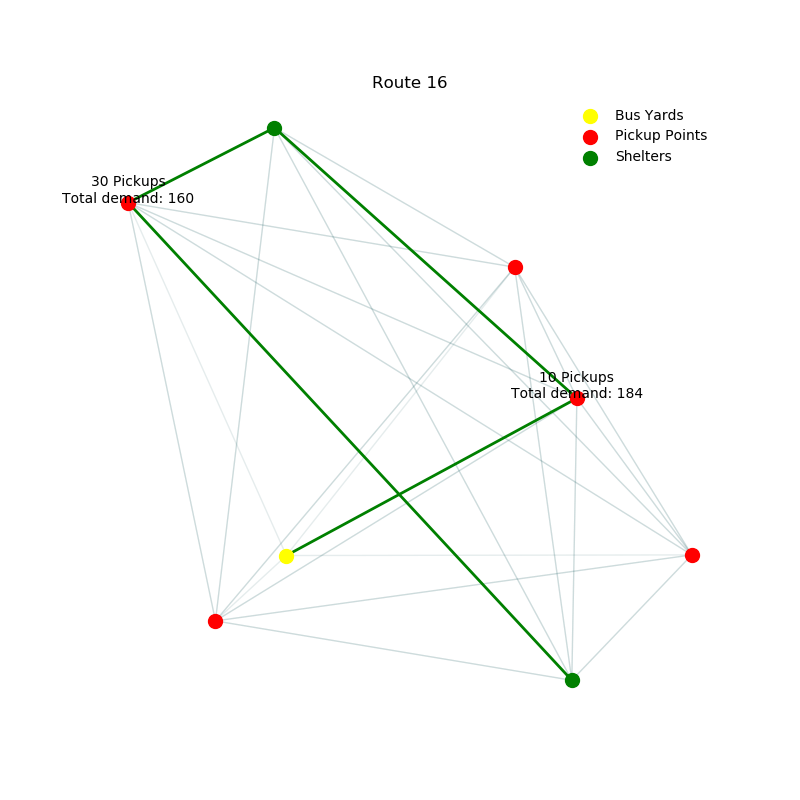}
\caption{Evacuation route of vehicle 17}
\label{P3}
\end{figure}

The simulation results of 7 instances is shown in Table II and compared in four categories, optimality gap, elapsed time, evacuation time, objective value. Note that 6 out of 7 instances obtain an optimal solution within prescribed elapsed time in which 5 out 7 instances attain its optimal solution under 150 seconds. The average optimality loss for 7 instances is $0.24\% $, and the highest optimality gap is $1.64\% $ occurred in Instance 8. The detailed evacuation routes for vehicle 1, vehicle 4, and vehicle 17 of instance 2  
are illustrated in Fig.\ref{P1}, Fig.\ref{P2}, and Fig.\ref{P3} respectively.

\par Comparing the simulation results obtained from  random and practical scenarios, we conclude that our model achieves good performance in solution optimality while the average loss of optimality in both scenarios is lower than $0.4\% $ within 3600 seconds. In terms of the average elapsed time, due to the difference in the structure of optimization problem, random scenario (1130.76s) is slightly higher than practical scenario (810.86s).

\section{Conclusion}

In this paper, to obtain a evacuation plan for evacuees, a bus based evacuation optimization problem is proposed and formulated as a mixed integer nonlinear programming (MINLP) problem which is a well-known NP-hard problem and computation time consuming. Then, unlike the cumbersome Big-M method, a Big-M free method is employed to simplify original MINLP problem by converting the original one as the mixed integer linear program (MILP) problem. Finally, the resultant MILP problem is solvable with efficient commercial solvers CPLEX. The effectiveness of proposed model is validated by extensive computational results.
\par In the future, more efficient algorithm, such as heuristics methods or decomposition methods, will be investigated to solve the proposed evacuation optimization problem.

\bibliographystyle{ieeetr}  
\bibliography{references.bib}






\end{document}